\documentclass[12pt]{amsart}
\usepackage{amsmath,amsthm,amsfonts,amssymb}
\begin{document} 
\newcommand{\Alb}{\mathop{\mathrm{Alb}}}
\newcommand{\m}{\mathfrak{m}}
\newcommand{\Spec}{{\mathrm{Spec}~}}
\newcommand{\Sing}{\mathop{\mathrm {Sing}}}
\newcommand{\supp}{\mathop{supp}}
\newcommand{\tensor}{\otimes}
\newcommand{\codim}{\mathop{codim}}
\newcommand{\B}{{\mathbb B}}
\newcommand{\C}{{\mathbb C}}
\newcommand{\N}{{\mathbb N}}
\renewcommand{\O}{{\mathcal O}}
\newcommand{\Q}{{\mathbb Q}}
\newcommand{\Z}{{\mathbb Z}}
\renewcommand{\P}{{\mathbb P}}
\newcommand{\R}{{\mathbb R}}
\newcommand{\rc}{\subset}
\newcommand{\rank}{\mathop{rank}}
\newcommand{\trace}{\mathop{tr}}
\newcommand{\dimc}{\mathop{dim}_{\C}}
\newcommand{\Lie}{\mathop{Lie}}
\newcommand{\Hom}{\mathop{\rm Hom}}
\newcommand{\Aut}{\mathop{{\rm Aut}}}
\newcommand{\Auto}{\mathop{{\rm Aut}_{\mathcal O}}}
\newcommand{\alg}[1]{{\mathbf #1}}
\newcommand{\res}{{\mathop{res}}}
\newtheorem*{definition}{Definition}
\newtheorem*{claim}{Claim}
\newtheorem{corollary}{Corollary}
\newtheorem*{Conjecture}{Conjecture}
\newtheorem*{GenAss}{General Assumptions}
\newtheorem{example}{Example}
\newtheorem{remark}{Remark}
\newtheorem*{observation}{Observation}
\newtheorem*{fact}{Fact}
\newtheorem*{remarks}{Remarks}
\newtheorem{lemma}{Lemma}
\newtheorem{proposition}{Proposition}
\newtheorem{theorem}{Theorem}
\newtheorem*{MainTheorem}{Main Theorem}
\title[Manifolds with trivial log tangent bundle]{%
On Manifolds with trivial logarithmic tangent bundle: The non-K\"ahler case.
}
%\shorttitle{}
\author {J\"org Winkelmann}
\begin{abstract}
We study non-K\"ahler manifolds with trivial logarithmic
tangent bundle. We show that each such manifold
arises as a fiber bundle with a compact complex parallelizable
manifold as basis and a toric variety as fiber. 
\end{abstract}
\address{%
J\"org Winkelmann \\
 Institut Elie Cartan (Math\'ematiques)\\
 Universit\'e Henri Poincar\'e Nancy 1\\
 B.P. 239, \\
 F-54506 Vand\oe uvre-les-Nancy Cedex,\\
 France
}
\email{jwinkel@member.ams.org\newline\indent{\itshape Webpage: }%
http://www.math.unibas.ch/\~{ }winkel/
}
\maketitle
\section{Introduction}
By a classical result of Wang \cite{W} a connected
compact complex manifold $X$
has holomorphically trivial tangent bundle if and only if there is
a connected complex Lie group $G$ and a discrete subgroup $\Gamma$ such that
$X$ is biholomorphic to the quotient manifold $G/\Gamma$.
In particular $X$ is homogeneous.
If $X$ is K\"ahler, $G$ must be commutative and the quotient manifold
$G/\Gamma$ is a compact complex torus.

Given a divisor $D$ on a compact complex manifold $\bar X$, we can define
the sheaf of {\em logarithmic differential forms}.
This is a coherent sheaf. If $D$ is locally s.n.c., this sheaf is
locally free and we may ask whether it is globally free,
i.e. isomorphic to $\O^{\dim\bar X}$.

In \cite{W2} we investigated this question in the case where the
manifold is K\"ahler (or at least in class ${\mathcal C}$) and
showed in particular that under these conditions there is a complex
semi-torus $T$ acting on $X$ with $X\setminus D$ as open orbit and
trivial isotropy at points in $X\setminus D$.

In this article we investigate the general case, i.e.,
we do not require the complex manifold $X$ to be K\"ahler or to be
in class ${\mathcal C}$.

There are two obvious classes of examples:
First, if $G$ is a connected complex Lie group with a discrete cocompact
subgroup $\Gamma$, then $X=G/\Gamma$ is such an example with $D$ being
the empty divisor.
Second, if $T$ is a semi-torus and $T\hookrightarrow\bar T$ is a smooth
equivariant compactification such that all the isotropy groups
are again semi-tori, then $X=\bar T$, $D=\bar T\setminus T$ yields
such examples (\cite{W2}).

Both classes contain many examples. E.g., every semisimple Lie group
admits a discrete cocompact subgroup (\cite{B}). On the other hand, every
semi-torus admits a smooth equivariant compactification such that
all isotropy groups are semi-tori. If the semi-torus under
consideration is algebraic, i.e.~a semi-abelian variety,
then this condition on the isotropy groups is fulfilled for
every smooth algebraic equivariant compactification.

We will show that in general every compact complex manifold
$\bar X$ with a divisor $D$ such that $\Omega^1(\log D)$ is trivial
must be constructed out of these two special classes.
First we show that there is a fibration $\bar X\to Y$
where the fiber is
an equivariant compactification $\bar A$
of a semi-torus $A$ and the base $Y$
is biholomorphic to a quotient of a complex Lie group $G$
by a discrete cocompact subgroup $\Gamma$.
From this one can deduce that any such manifold $X$ is a compactification
of a $(\C^*)^l$-principal bundle over a compact complex parallelizable
manifold, just like semi-tori are extension of compact complex tori
by some $(\C^*)^l$.

\section{The main results}

\begin{theorem}\label{thm1}
Let $\bar X$ be a compact complex manifold and let $D$ be a 
locally s.n.c. divisor such that $\Omega^1(\log D)$ is globally
trivial (i.e.~isomorphic to ${\mathcal O}_{\bar X}^{\dim\bar X}$).

Let $G$ denote the connected component of
\[
\Aut(\bar X,D)=\{g\in\Aut(\bar X):g(D)=D\},
\]
and let $x\in X=\bar X\setminus D$.

Then $G$ acts transitively on $X=\bar X\setminus D$ with discrete
isotropy group $\Gamma=\{g\in G:g\cdot x=x\}$.

Furthermore there exists a central subgroup
$C\simeq(\C^*)^l$ of $G$, a smooth equivariant compactification
$C\hookrightarrow \bar C$,
a compact complex parallelizable manifold $Y$
and a  holomorphic fibration $\pi:\bar X\to Y$
such that

\begin{itemize}
\item
$\pi$ is a locally holomorphically trivial fiber bundle
with $\bar C$ as fiber and $C$ as structure group.
\item
The $\pi$-fibers are closures of orbits of $C$,
\item
The projection map $\pi$ is $G$-equivariant and admits
a $G$-invariant flat holomorphic connections.
\item
There is a connected complex Lie subgroup $H\subset G$ 
acting transitively on $Y$ with discrete isotropy group
$\Lambda$.
\item
$\Lie(G)=\Lie(H)\oplus\Lie(C)$ and
$\Gamma=\{\gamma\cdot\rho(\gamma):\gamma\in\Lambda\}$
for some group homomorphism $\rho:\Lambda\to C$.
\end{itemize}
\end{theorem}

Conversely, any such fiber bundle yields a log-parallelizable
manifold:

\begin{theorem}\label{thm2}
Let $Y$ be a compact complex parallelizable manifold,
$C\simeq(\C^*)^l$, $\bar C$ a smooth equivariant compactification
of $C$. 
Assume that all the isotropy groups of the $C$-action
on $\bar C$ are reductive.
Let $E_1,\ldots, E_l$ be $\C^*$-principal bundles over $Y$,
each endowed with a holomorphic connection.

Define $\bar X$ as the total space of the $\bar C$-bundle associated
to the $C$-principal bundle $\Pi E_i$ and let $D$ be the divisor
on $\bar X$ induced by the divisor $\bar C\setminus C\subset\bar C$.

Then $D$ is a locally s.n.c.~divisor and
$\Omega^1(\bar X; \log D)$ is globally trivial.
\end{theorem}

This can also be expressed in group-theoretic terms:
\begin{theorem}\label{thm3}
Let $H$ be a complex connected Lie group, $\Lambda$ a discrete
cocompact subgroup, $l\in\N$, $C=(\C^*)^l$, $i:C\hookrightarrow
\bar C$ a smooth equivariant compactification and $\rho:\Lambda
\to C$ a group homomorphism.
Assume that all the isotropy groups of the $C$-action
on $\bar C$ are reductive.

Define $\bar X=\left( H\times\bar C\right)/\!\sim$ and
$X=\left(H\times C\right)/\!\sim$ where
$(h,x)\sim (h',x')$ iff there is an element $\lambda\in\Lambda$
such that 
$\left(h',x'\right)=\left(h\lambda,\rho(\lambda^{-1})x\right)$.
Define $D=\bar X\setminus X$.

Then $D$ is a locally s.n.c.~divisor and
$\Omega^1(\bar X;\log D)$ is globally trivial.
\end{theorem}

\section{Logarithmic forms}
Let $\bar X$ be a complex manifold and $D$ a divisor.
We say that $D$ is ``locally s.n.c'' (where s.n.c.~stands
for ``simple normal crossings'') if
if for every point $x\in \bar X$ there exists local
coordinates $z_1,\ldots,z_n$ and a number $d\in\{0,\ldots n\}$ such that
in a neighbourhood of $x$ the divisor $D$ equals the zero divisor
of the holomorphic function $\Pi_{i=1}^d z_i$.

It is called a ``divisor with only simple normal crossings as 
singularities'' or ``s.n.c.~divisor'' if in addition every irreducible
component of $D$ is smooth.

Let $\bar X$ be a compact complex manifold with 
a locally s.n.c.{} divisor $D$.
There is a stratification as follows: $X_0=X=\bar X\setminus D$,
$X_1=D\setminus \Sing(D)$ and for $k>1$ the stratum $X_k$ is the
non-singular part of $\Sing(\bar X_{k-1})$.
If in local coordinates $D$ can be written as $\{z:\Pi_{i=1}^d z_i=0\}$,
then $z=(z_1,\ldots,z_n)\in X_k$ iff $\#\{i:1\le i \le d, z_i=0\}=k$.

Let $D$ be an effective divisor on a complex manifold $\bar X$.
The sheaf $\Omega^k(\log D)$ of 
{\em logarithmic $k$-forms} with respect to $D$ is defined as
the sheaf of all meromorphic $k$-forms $\omega$ for which
both $f\omega$ and $fd\omega$ are holomorphic for all $f\in
{\mathcal I}_D$.
If $D$ is locally s.n.c., we may also define the algebra
(with respect to exterior product) of logarithmic $k$-forms
as the $\O_{\bar X}$-algebra generated by $1$ and all
$df/f$ where $f$ is a section ${\mathcal O}_{\bar X}\cap{\mathcal O}_X^*$.

Note: A logarithmic $0$-form is simply a holomorphic $0$-form, i.e.,
a holomorphic function.

For every natural number $k$ the
sheaf of logarithmic $k$-forms is
a coherent ${\mathcal O}$-module sheaf. It is locally free if
$D$ is a locally s.n.c.{} divisor. Especially, if
$D=\{z_1\cdot\ldots\cdot z_d=0\}$, then $\Omega^1(\log D)$
is locally the free $\O_{\bar X}$-module over 
$dz_1/z_1,\ldots,dz_d/z_d,dz_{d+1},\ldots,dz_n$.

For a {\em locally s.n.c.} divisor $D$ on a complex manifold
we define the {\em logarithmic
tangent bundle} $T(-\log D)$ as the dual bundle of $\Omega^1(\log D)$.

Then $T(-\log D)$ can be identified
with the sheaf of those holomorphic vector fields $V$ on $\bar X$ 
which fulfill the following property:

{\em $V_x$ is tangent to $X_k$ at $x$ for every $k$ and 
every $x\in X_k$}.

In local coordinates: If $D=\{z:\Pi_{i=1}^d z_i=0\}$, then
$T(-\log D)$ is the locally free sheaf generated by
the vector fields $z_i\frac{\partial}{\partial z_i}$ ($1\le i\le d$)
and $\frac{\partial}{\partial z_i}$ ($d<i\le n$).

This implies in particular (compare also \cite{Wi}):
\begin{proposition}\label{prop-almost-homo}
Let $X$ be a compact complex manifold and $D$ a locally s.n.c. divisor
such that $\Omega^1(\log D)$ is globally trivial.

Let $G$ denote the connected component of the group
$\Aut(X,D)$ of all holomorphic automorphisms of $X$ which stabilize
$D$.

Then the Lie algebra of $G$ can be identified with
$\Gamma(\bar X,\Omega^1(\log D))^*$ and the $G$-orbits in $\bar X$
coincide with the connected components of the strata $X_k$.

For every $x\in X_0=\bar X\setminus D$ the isotropy group
$G_x=\{g\in G:g(x)=x\}$ is discrete.

The isotropy representation of $G_x$ on $T_x$ is almost faithful
(i.e.~the kernel is discrete) for every $x\in\bar X$.
\end{proposition}
\section{Residues}
Let $\bar X$ be a complex manifold and $D$ a locally s.n.c. divisor.
In the preceding section we introduced the notion of logarithmic
forms. Now we introduce ``residues''.

Let $\tau:\hat D\to D$ denote a normalization of $D$.
Note that $\hat D$ is a (not necessarily connected) compact manifold.
It is smooth, because $D$ is locally s.n.c. Furthermore
$D$ being locally s.n.c.~implies that
$D'=\tau^*\Sing(D)$ defines a divisor on $\hat D$.

We are going to define 
\[
res: \Omega^k(X;\log D) \longrightarrow \Omega^{k-1}(\hat D,\log D').
\]

Let $\hat x\in \hat D$ and $x=\tau(\hat x)\in \bar X$.
$x$ admits an open neighbourhood $W$ in $\bar X$ with local
coordinates $z_1,\ldots,z_n$ such that
\begin{enumerate}
\item $z_1(x)=\ldots=z_n(x)=0$.
\item
$\tau$ induces an isomorphism between an open neighbourhood $\hat W$
of $\hat x$ in $\hat D$ with $\{p\in W:z_1(p)=0\}$.
\item
There is a number $1\le k \le n$ such that $D|_W$ is the zero divisor
of the function $z_1\cdot\ldots\cdot z_k$.
\item
The restriction of $D'$ to $\hat W$ is the pull-back of the zero
divisor of the function $z_2\cdot z_3\cdot\ldots\cdot z_k$ on $W$.
\end{enumerate}

Now let $\omega\in\Gamma(W,\Omega^k(\log D))$.
Then 
\[
\omega=\eta+\frac{dz_1}{z_1}\wedge\mu
\]
where $\eta$ and $\mu$ have poles only along the zero set
of $z_2\cdot\ldots\cdot z_k$.
We define
the restriction of
$\res(\omega)$ to $\hat W$ as
\[
\res(\omega)|_{\hat W} = \tau^*\mu
\]
It is easily verified that $\res(\omega)$ is independent of the choice
of local coordinates and the decomposition 
$\omega=\eta+\frac{dz_1}{z_1}\wedge\mu$ and that furthermore
the only poles of $\res(\omega)$ are logarithmic poles along $D'$.

Moreover:
\begin{lemma}
Let $X$ be a complex manifold, $D$ a locally s.n.c. divisor
and $\omega$ a logarithmic $k$-form ($k\in\N\cup\{0\}$).

Then
$\omega$ is holomorphic iff $\res(\omega)=0$.
Furthermore 
\[
res(d\omega) = - d\,\res(\omega)
\]
\end{lemma}
\begin{proof}
This is easily verified by
calculations in local coordinates.
\end{proof}
As a consequence we obtain the following result which we will
need later on.
\begin{lemma}\label{d-no-poles}
Let $\bar X$ be a compact complex manifold, $D$ a locally s.n.c.
divisor and $\omega$ a logarithmic $1$-form.

Then $d\omega$ is holomorphic.
\end{lemma}
\begin{proof}
The residue $\res(\omega)$
is a logarithmic $0$-form on the normalization $\hat D$ of $D$.
But logarithmic $0$-forms are simply holomorphic functions, and
every holomorphic function on a $\hat D$ is locally constant
by our compactness assumption. Therefore $d\res(\omega)=0$.
Consequently $\res(d\omega)=0$,
which is equivalent to the condition that $d\omega$ has no poles.
\end{proof}
\begin{remark}
Residues can also be introduced for divisors which are not locally
s.n.c., see (Saito?).
\end{remark}

\section{Group theoretic preparations}
\begin{definition}
A (not necessarily closed) subgroup $H$ of a topological group $G$
is called {\em cocompact} iff there exists a compact subset $B\subset G$
such that $B\cdot H=\{bh:b\in B,h\in H\}=G$.
\end{definition}
If $H$ is a closed subgroup of a Lie group $G$, then $H$ is cocompact if and 
only if the quotient space $G/H$ is compact.

\begin{proposition}
Let $G$ be a connected complex Lie group, $Z$ its center and $\Gamma$
a subgroup. Let $C=Z_G(\Gamma)^0$ denote the connected component
of the centralizer
$Z_G(\Gamma)=\{g\in G;g\gamma=\gamma g\ \forall\gamma\in\Gamma\}$.

Assume that $Z\Gamma$ is cocompact.

Then $C=Z^0$.
\end{proposition}
\begin{proof}
Let $B$ be a compact subset of $G$ such that $BZ\Gamma=G$.

For each $c\in C$ we define a holomorphic map from $G$ to $G$ 
by taking the commutator with $c$:
\[
\zeta_c:g\mapsto cgc^{-1}g^{-1}
\]

Since $c$ commutes with every element of $Z\Gamma$, we have
$\zeta_c(G)=\zeta_c(B)$.
Thus compactness of $B$ implies that
the image $\zeta_c(G)=\zeta_c(B)$ is compact.
Considering the adjoint representation
\[
Ad:G\to GL(\Lie G)\subset\C^{n\times n}
\]
we deduce that for every holomorphic function $f$ on $GL(\Lie G)$
the composed map $f\circ Ad\circ\zeta_c:G\to\C$ is a bounded
holomorphic function on $G$ and therefore constant.
It follows that the image $\zeta_c(G)$ is contained in $\ker Ad=Z$.

%Now $Z^0$ is a connected commutative complex Lie group. Let $U$ denote
%ist maximal (real) compact Lie subgroup.
%For every $x\in Z^0\setminus U$ there is a real Lie group
%homomorphism $\rho:Z^0\to\R$ with $\rho(x)\ne 0$.
%Now $\rho\circ\exp:\Lie Z\to\R$ is just a $\R$-linear map.
%It follows that $\rho$ and therefore $\rho\circ\zeta_c$ is
%a plurisubharmonic function on $G$ for every $c\in C$ and 
%$\rho\in\Hom(Z^0,\R)$. Using the compactness of $\zeta_c(G)$
%it follows that $\rho\circ\zeta_c$ is a bounded plurisubharmonic
%function and therefore constant.
%Hence $\zeta_c(G)\subset U$.
Thus $\zeta_c(h)$ is central in $G$ for all $c\in C$, $h\in G$.
This implies
\begin{multline*}
\zeta_c(g)\zeta_c(h^{-1})
= c g c^{-1} g^{-1} (c h^{-1} c^{-1} h)
= c g c^{-1} (c h^{-1} c^{-1} h) g^{-1} \\
= c g h^{-1} c^{-1} h g^{-1}
=\zeta_c(gh^{-1})
\end{multline*}

for all $c\in C$ and $g,h\in G$.
Hence $\zeta_c:G\to G$ is a group homomorphism for all $c\in C$.

It follows that $\zeta_c(G)$ is a connected
compact complex Lie subgroup of $Z$.
In particular, $\zeta_c(G)$ is a compact complex torus
and the commutator group $G'$ is contained in $\ker\zeta_c$.
Let $A$ denote the smallest closed complex Lie subgroup of $G$ containing
$G'Z\Gamma$. Then $A$ is normal, because it is a subgroup containing
$G'$, and $G/A$ is compact, since $Z\Gamma$ is cocompact.
Therefore $G/A$ is a connected compact complex Lie group, i.e., a compact
complex torus.
Thus for every $c$ the map $\zeta_c:G\to G$ fibers through $G/A$
and is induced by some holomorphic Lie group homomorphism
from $G/A$ to $T$ where $T$ denotes the maximal compact complex
subgroup of $Z$.
Note that both $G/A$ and $T$ are compact tori.
This implies that $\Hom(G/A,T)$ is discrete.
Now
$c\mapsto\zeta_c$ defines a map from the {\em connected} complex Lie group
$C$ to $\Hom(G/A,T)$.
Since
$\Hom(G/A,T)$ is discrete, this map is constant.
Furthermore $\zeta_e\equiv e$.
Hence $\zeta_c\equiv e$ for all $c\in C$.
This implies $C\subset Z$. On the other hand the inclusion $Z^0\subset C$
is obvious. This completes the proof of the equality $C=Z^0$.
\end{proof}

\begin{proposition}\label{prop-z-closed}
Let $G$ be a connected complex Lie group, $Z$ its center and $\Gamma$
a discrete subgroup. 
Assume that $Z\Gamma$ is cocompact.

Then $Z^0\Gamma$ is closed in $G$.
\end{proposition}
\begin{proof}
Let $C=C_G(\Gamma)^0$ be the connected component of the
centralizer. Then $C=Z^0$ be the preceding proposition.
On the other hand $C\Gamma$ is closed by an argument of
Raghunathan, see \cite{Wi}, lemma 3.2.1.
\end{proof}

\begin{corollary}\label{cor-ext}
Let $G$ be a connected complex Lie group, $Z$ its center and $\Gamma$
a discrete subgroup. 
Assume that $Z\Gamma$ is cocompact.

Then there exists a discrete cocompact subgroup $\Gamma_1$ in $G$
with $\Gamma\subset\Gamma_1$. Moreover the commutator groups
$\Gamma'$, $\Gamma_1'$ can be required to be equal.
\end{corollary}
\begin{proof}
Note that $Z^0$ is a connected commutative Lie group. It is easy to see
that there is a discrete subgroup $\Lambda\subset Z^0$ which is 
cocompact in $Z^0$ and contains $\Gamma\cap Z^0$.
Define $\Gamma_1=\Gamma\cdot\Lambda$.
This is a subgroup, because $\Lambda$ is central in $G$.
By prop.~\ref{prop-z-closed} both $Z^0\Gamma_1$ and $Z^0\Gamma$ are closed
in $G$.
Hence we may consider the fibration
\[
X=G/\Gamma_1 \to G/Z^0\Gamma_1=G/Z^0\Gamma=Y.
\]
Now $Y$ and $Z^0/\Lambda\simeq Z^0\Gamma_1/\Gamma_1$  are both compact,
hence $X$ is likewise compact.

Finally we note that the equality $\Gamma'=\Gamma_1'$ follows from the
fact that $\Lambda$ is central.
\end{proof}

\begin{proposition}
Let $G$ be a connected complex Lie group, $Z$ its center and $\Gamma$
a discrete subgroup. 
Assume that $Z\Gamma$ is cocompact.

Then $\Gamma'\cap Z$ is cocompact in $G'\cap Z$.
\end{proposition}
\begin{proof}
By cor.~\ref{cor-ext} there is a discrete cocompact subgroup
$\Gamma\subset\Gamma_1\subset G$ with $\Gamma'=\Gamma_1'$.

By passing to the universal covering there is no loss in generality
in assuming that $G$ is simply-connected.

Now $(G'\cap R)/(\Gamma_1\cap G'\cap R)$ and $Z/(Z\cap\Gamma_1)$
are both compact, hence $(G'\cap Z)/(G'\cap Z\cap\Gamma_1)$
is compact, too.
Recall that $\Gamma_1'=\Gamma'\subset\Gamma$. By
\cite{Wi}, prop.~3.11.2 
$R\cap\Gamma'_1$ is of finite index in $G'\cap R\cap\Gamma_1$.
Using $Z\subset R$, it follows that $Z\cap G'\cap\Gamma$ is of
finite index in $Z\cap G'\cap\Gamma_1$. Hence
$Z\cap G'\cap\Gamma$ is cocompact in $Z\cap G'$.
\end{proof}
\begin{corollary}\label{cor-der-cpt}
Let $G$ be a connected complex Lie group, $Z$ its center and $\Gamma$
a discrete subgroup. Assume that $Z\Gamma$ is cocompact and assume
further that $G$ acts effectively on $G/\Gamma$.

Then $G'\cap Z$ is compact.
\end{corollary}
\begin{proof}
The subgroup $\Gamma\cap Z$ acts trivially on $G/\Gamma$, hence
$\Gamma'\cap Z=\{e\}$ if $G$ acts effectively on $G/\Gamma$.
However, the trivial subgroup $\{e\}$ is cocompact in $G'\cap Z$
if and only if the latter is compact.
\end{proof}
\section{The commutative case}

\begin{proposition}\label{commutative}
Let $G$ be a commutative connected complex Lie group
and let $G\hookrightarrow \bar X$ be a smooth equivariant compactification,
$D=\bar X\setminus G$. 

Then $\Omega^1(\log D)$ is a holomorphically trivial vector bundle
on $\bar X$ if and only if the following two conditions are fulfilled:
\begin{enumerate}
\item $G$ is a semi-torus, and
\item every isotropy group for the $G$-action on $\bar X$ is a semi-torus.
\end{enumerate}
\end{proposition}
\begin{proof}
This is $(1)\Rightarrow(2)$ of the Main Theorem in \cite{W2}.
In \cite{W2} we assumed that $\bar X$ is K\"ahler or in class ${\mathcal C}$.
However, this assumption is used only in order to {\em deduce} that $G$
is commutative. Hence the line of arguments in \cite{W2} still applies
if instead of requiring $\bar X$ to be K\"ahler we merely
assume that $G$ is 
commutative.
\end{proof}

\section{The isotropy groups}
\begin{proposition}\label{iso-central}
Let $\bar X$ be a compact complex manifold
and $D$ a locally s.n.c. divisor such that $\Omega^1(\bar X;\log D)$
is spanned by global sections.
Let $G$ be a connected complex Lie group acting holomorphically
on $\bar X$ such that $D$ is stabilized and $p\in\bar X$.

Then the connected component of the isotropy group 
$G_p=\{g\in G:g(p)=p\}$ is central in $G$.
\end{proposition}
\begin{proof}
Let $w\in\Lie(G)$ and $v\in\Lie(G_p)$.
We have to show that $[w,v]=0$ for any such $w,v$.
By abuse of notation we identify the elements $w,v$ of $\Lie(G)$
with the corresponding fundamental vector fields on $\bar X$.
Now let us assume that $[w,v]\ne 0$. Then there exists a logarithmic
one-form $\omega\in\Omega^1(\bar X;\log D)$ such that
$(\omega,[w,v])$ is not identically zero.
Since $G$ stabilizes $D$, the vector fields $w,v$ are tangent to $D$ and
thus sections in $T(-log D)$ which is dual to $\Omega^1(\log D)$.
Therefore $(\omega,w)$ and $(\omega,v)$ are global holomorphic functions
on $\bar X$. Because $\bar X$ is compact, they are constant.
Hence $w(\omega,v)$ and $v(\omega,w)$ vanish both.
Therefore
\[
\omega([w,v]) = d\omega(w,v).
\]
Now $d\omega$ has no poles (lemma~\ref{d-no-poles}) 
and $v$ vanishes at one point,
namely $p$. Thus $\omega([w,v])$ is a global holomorphic function on
a compact connected manifold which vanishes at one point, i.e.,
$\omega([w,v])\equiv 0$ contrary to our assumption on $\omega$.
We have thus deduced that the assumption $[w,v]\ne 0$ leads to a
contradiction. Hence $[w,v]=0$ for all $w\in\Lie(G)$ and all
$v\in\Lie(G_p)$. It follows that $G_p^0$ is central in $G$.
\end{proof}

\section{Cocompactness of $Z\Gamma$}
We fix some notation to be used throughout this section:
\begin{itemize}
\item
$\bar X$ is a compact complex manifold,
\item $D$ a loc.~s.n.c.~diuvisor on $\bar X$ and $X=\bar X\setminus|D|$,
\item We assume that $\Omega^1(\bar X,\log D$ is a globally
trivial ${\mathcal O}_{\bar X}$-module.
\item $\Aut(\bar X,D)=\{g\in\Aut(\bar X):g(D)=D\}$,
\item $G$ denotes the connected component of $\Aut(\bar X,D)$ which contains
the identity map.
\item $x$ denotes a fixed point in $X=\bar X\setminus|D|$,
\item $\Gamma$ denotes the isotropy group at $x$, i.e.,
$\Gamma=\{g\in G:g\cdot x=x\}$
\item $Z$ denotes the center of $G$.
\end{itemize}

\begin{proposition}\label{zcocompact}
$Z\Gamma$ is cocompact in $G$.
\end{proposition}
\begin{proof}
Equip $\bar X$ with some Riemannian metric.
Let $p\in |D|$. Assume that $p\in D_1\cap\ldots D_j$ and $p\not\in D_i$ for
$i>j$. Then there are logarithmic vector fields $v_i$ and local
holomorphic cooordinates $z_i$ near $p$ such that 
$v_i=z_i\frac{\partial}{\partial z_i}$ for $1\le i\le j$.
As a consequence, for every point $p\in \bar X$ there exists an open
neighborhood $W_p$ and a positive number $\epsilon_p>0$
such that the following assertion holds:
{\em For every $y\in W_p\setminus |D|$ 
there exists an element $g\in G_p^0$
such that $d(g\cdot w,|D|)>\epsilon_p$.}
Since $\bar X$ is compact, it can be covered by finitely many of the
open sets $W_p$. We choose such a finite collection and define $\epsilon$
as the minimum of the $\epsilon_p$.
Let $H$ be the subgroup of $G$ generated by all the subgroups $G_p^0$.

Then we have: {\em For every $y\in X\setminus |D|$ there exists an element
$g\in H$ such that $d(g\cdot y,|D|)>\epsilon$.}

Compactness of $X$ implies that $C=\{z\in X:d(z,|D|)\ge\epsilon\}$
is compact. We have $H\cdot C=X\setminus|D|$. 
We can choose a compact subset $K\subset G$ such that
the natural projection from $G$ onto $G/\Gamma\simeq X\setminus|D|$
maps $K$ onto $C$. Then $G=H\cdot K\cdot \Gamma$.
Observe that $H$ is
central in $G$ due to prop~\ref{iso-central}.
Hence 
\[
G=H\cdot K\cdot \Gamma=K\cdot H\cdot \Gamma
\]
and $H\subset Z$.
It follows that $Z\Gamma$ is cocompact in $G$.
\end{proof}
\begin{corollary}\label{cor-z}
$Z^0\Gamma$ is closed and cocompact in $G$.
\end{corollary}
\begin{proof}
Implied by prop.~\ref{prop-z-closed} and \ref{zcocompact}.
\end{proof}

\begin{proposition}\label{prop-fibration}
$G/Z^0\Gamma$ is a compact parallelizable manifold
and the natural projection $G/\Gamma\to G/Z^0\Gamma$
extends to a holomorphic map from $\bar X$ onto $G/Z^0\Gamma$.
\end{proposition}
\begin{proof}
$G/Z^0\Gamma$ can be described as the quotient of the complex
Lie group $G/Z^0$ by the discrete subgroup $\Gamma/(\Gamma\cap Z^0)$
and is therefore a compact parallelizable manifold.

Now let $p$ be a point in $D\setminus\Sing(D)$.
Choose local cooordinates near $p$ such that $D=\{z_1=0\}$.
Since $\Omega^1(\log D)$ is trivial, there is a global vector field
$v$ such that 
\[
(\frac{1}{z_1}v)_p={\frac{\partial}{\partial z_1}}_p.
\]
Now $v$ is central in $\Lie(G)$ and $G$ acts locally transitively on
$D\setminus\Sing(D)$. Therefore $v$ vanishes as holomorphic vector field
on $D$. It follows that $w=\frac{1}{z_1}v$ is a holomorphic vector field
near $p$.
Moreover $w_p={\frac{\partial}{\partial z_1}}_p$.
Now $w$, as any vector field, is locally integrable.
Therefore there is a local coordinate system near $p$ in which
$w$ equals $\frac{\partial}{\partial z_1}$.
Since $v$ is central as element in $\Lie(G)$, $v$, as a holomorphic
vector field on $X\setminus D$, is tangent to the fibers of
$\tau:G/\Gamma\to G/Z^0\Gamma$. Thus in this local coordinate
system near $p$, the map $\tau$ depends only on the variables 
$z_2,\ldots, z_n$. But this implies that $\tau$, which is defined on
$\{z_1\ne 0\}$, extends to a map defined on the whole neighbourhood
of $p\in D$. In this way we say that $\tau$, originally defined on
$X\setminus D$, extends to a holomorphic map defined on 
$X\setminus\Sing(D)$.

Next we observe that the image space $G/Z^0\Gamma$ is complex parallizable,
and therefore has a Stein universal covering. It follows that we can use
the classical Hartogs theorem to extend $\tau$ through the set $\Sing(D)$
which has codimension at least two in $X$.
\end{proof}

\begin{proposition}\label{fib-properties}
Let $Y$ be a complex manifold on which $G$ acts transitively,
and $\pi:X\to Y$ an equivariant holomorphic map.

Then each fiber $F$ of $\pi$ is smooth, $D\cap F$ is a s.n.c. divisor
in $F$ and $\Omega^1(F,\log (D\cap F))$ is trivial.
\end{proposition}
\begin{proof}
A generic fiber is smooth. Since $\pi$ is equivariant and $G$ acts
transitively on $Y$, it follows that every fiber is smooth.
Let $S=\{y\in Y:\pi^{-1}(y)\subset D\}$. Then $S$ is $G$-invariant.
Hence either $\pi(S)=Y$ or $S=\emptyset$. But $\pi(S)=Y$ would imply
$D=X$ contrary to $D$ being a divisor. Hence $S=\emptyset$, i.e.
$D$ does not contain any fiber of $\pi:X\to Y$.

Now we consider the usual stratification of $D$ (as described in \S3).
At each point $p\in D_k$ the divisor $D$ can locally be defined
as $D=\{f_1\ldots f_k=0\}$ where $df_1,\ldots,df_k$ are linearly
independent in $T^*_pX$. Since $D$ can not be a pull-back of a
divisor on $Y$, for a generic choice of $p$ there will be local functions
$f_1,\ldots,f_k$ near $p$ such that $D=\{f_1\cdot\ldots\cdot f_k=0\}$
and in addition such that the $df_i$ are linearly independent as elements
in $T_pF$ (with $F=\pi^{-1}(\pi(p))$). In the same spirit as in the definition
of $S$ above we can now define the set of points in $D$ where this fails.
Observing that these sets are invariant, and recalling that
$G$ acts transitively on the connected components of the
strata $D_k$, we deduce that they must 
empty.
Thus $D\cap F$ is always a s.n.c. divisor.

By similar arguments we see that $\Omega^1(F,\log(D\cap F))$
is trivial.
\end{proof}

\section{Proofs for the theorems}
\begin{proof}[Proof of theorem \ref{thm1}]
Let $Z$ denote the center of $G$.
By definition $G$ is a subgroup of the automorphism group
of $\bar X$ and therefore acts effectively on $\bar X$.
By cor.~\ref{cor-z} $Z^0\Gamma$ is closed in $G$ and
due to prop.~\ref{prop-fibration} there is a holomorphic fiber bundle
$\pi:\bar X\to Y_0=G/Z^0\Gamma$ which extends the natural projection 
$G/\Gamma\to Y_0=G/Z^0\Gamma$.
 From prop~\ref{fib-properties} we deduce that a typical fiber $F$ of $\bar\pi$
is log-parallelizable with respect to the divisor $F\cap D$.
Since $Z$ is commutative, this implies (due to prop.~\ref{commutative}) 
that $Z^0$ is a semi-torus.
Let $C$ be a maximal connected linear subgroup of $Z^0$.
Then $C\simeq(\C^*)^l$ for some $l\in\N$ and $Z^0/C$ is a compact
complex torus. 
The fibration $Z^0\mapsto Z^0/C$ induces a tower of
fibrations
\[
G/\Gamma \to G/C\Gamma \to G/Z^0\Gamma
\]
which (using prop.~\ref{prop-fibration}) extends to fibrations
\[
\bar X \to Y=G/C\Gamma \to G/Z^0\Gamma
\]
Now $Y$ is parallelizable because $C$ is normal in $G$ and compact
because both $G/Z^0\Gamma$ and $Z^0/C$ are compact.

From cor.~\ref{cor-der-cpt} we deduce that $G'\cap Z$ is compact.
Since $C\subset Z$ and $C\simeq(\C^*)^l$, it follows that $G'\cap C$
is discrete.
Thus we can choose a complex vector subspace $V\subset\Lie(G)$
such that $\Lie(G')\subset V$ and $\Lie(G)=V\oplus\Lie(C)$.
The condition $\Lie(G')\subset V$ implies that $V$ is an ideal
in $\Lie(G)$. Hence $V$ is the Lie algebra of a normal Lie subgroup $H$
of $G$. Furthermore $H\cdot C=G$ due to $\Lie(G)=V\oplus\Lie(C)$.
The condition $H\cdot C=G$ implies that $H$ acts transitively on
$Y=G/C\Gamma$. For dimension reasons the isotropy group
$\Lambda=H\cap(C\Gamma)$ is discrete.

Because $G$ acts effectively on $X$, the intersection $Z\cap\Gamma$
must be trivial. Hence $C\cap\Gamma=\{e\}$. It follows that the
projection map $\tau:G\to G/C\simeq H/(H\cap C)$ maps $\Gamma$
injectively onto $\tau(\Gamma)=C\Gamma/C$.

We claim that $\tau(\Gamma)=\tau(\Lambda)$. Indeed, 
$\Lambda=H\cap(C\Gamma)\subset C\Gamma$ implies 
$\tau(\Lambda)\subset\tau(\Gamma)$. On the other hand,
if $\gamma\in\Gamma$, then there exists an element $c\in C$
such that $\gamma c\in H$, because $H\cdot C=G$.
Then $\tau(\gamma)=\tau(\gamma c)$ and $\gamma c=c\gamma\in H\cap (C\Gamma)
=\Lambda$. Hence $\tau(\Gamma)\subset\tau(\Lambda)$.

It follows that there exists a map $\rho:\Lambda\to C$ such that
for each $\lambda\in\Lambda$ the product $\lambda\rho(\lambda)$
is the unique element $\gamma$ of $\Gamma$ with $\tau(\lambda)=\tau(\gamma)$.

One verifies easily that $\rho$ is a group homomorphism, using the
facts that $\tau$ is a group homomorphism and that $C$ is central.

To obtain the $G$-invariant flat connection, we observe that
the fundamental vector fields of $\Lie(H)=V$ induce a decomposiition
$T_x\bar X=V\oplus\ker(d\pi)$ in each point $x\in\bar X$.
This connection is obviously $G$-invariant. Moreover it is flat,
because $V=\Lie(H)$ is a {\em Lie subalgebra} of $\Lie(G)$.

This completes the proof.
\end{proof}

\begin{proof}[Proof of theorem~\ref{thm2}]
Due to prop.~\ref{commutative} we know that $\bar C$ is
log-parallelizable, i.e., that $\Omega^1(\bar C,\log (\bar C\setminus C))$
is isomorphic to ${\mathcal O}_{\bar C}^l$.
Moreover, this isomorphism is given by meromorphic $1$-forms
which are dual to a basis of the $C$-fundamental vector fields.
Let $V_1,\ldots,V_l$ be such a basis of $C$-fundamental vector fields
and $\eta_1,\ldots,\eta_l$ a dual basis for the logarithmic one-forms
on $\bar C$. We may regard $\bar C$ as one fiber of the projection map
from $\bar X$ onto $Y$. Using the $C$-principal right action of
$C$ we extend $V_1,\ldots,V_l$ to holomorphic $C$-fundamental vector
fields on all of $\bar X$. Let $H\subset T\bar X$ be the horizontal
subbundle defined by the connection.
Then we can extend the meromorphic one-forms $\eta_i$ as follows:
For each $i$ we require that $\eta_i(Y_j)=\delta_{ij}$ and that
$\eta_i$ vanishes on $H$.

Let $\omega_1,\omega_r$ be a family of holomorphic $1$-forms
on $Y$ which gives a trivialization of the tangent bundle $TY$.
Define $\mu_i=\pi^*\omega_i\in\Omega^1(\bar X)$.

We claim that the family $(\mu_i)_i$ together with the family $\eta_i$
gives a trivialization of the sheaf of logarithmic one-forms on $\bar X$.

Evidently the $\mu_i$ are holomorphic. Since the restriction of
$\eta_i$ to a fiber is logarithmic and $\eta_i|_H\equiv 0$,
it is clear that $f\eta_i$ is holomorphic for any locally given function
$f$ vanishing on $D$. It remains to show that $fd\eta_i$ is
holomorphic as well. To see this, we calculate $fd\eta_i(Y,Z)$
for holomorphic vector fields $Y,Z$.
Recall that
\[
d\eta_i(Y,Z)=Y\eta_iZ-Z\eta_iY-\eta_i([Y,Z])
\]
It suffices to verify holomorphicity for a base of vector fields.
Thus we may assume that each of the vector fields $Y$ and $Z$ is
horizontal or vertical (with respect to the connection).
If both are vertical, there is no problem since $\eta_i$ restricted
to a fiber is logarithmic.
If both are horizontal, then $\eta_iY=\eta_iZ=0$. Furthermore
$[Y,Z]$ is horizontal, because the connection is flat.
Hence $\eta_i([Y,Z])=0$ and consequently $d\eta_i(Y,Z)=0$.
Finally let us discuss the case where $Y$ is horizontal 
and $Z$ is vertical. 
Then $\eta_iY=0$ and moreover $[Y,Z]=0$
because $H$ is defined by a connection for the $C$-bundle and therefore
$C$-invariant. Thus $fd\eta_i(Y,Z)=fY\eta_iZ$.
Since the fiber $\bar C$ is log-parallelizable,
it suffices to consider the case where, up to multiplication by a
meromorphic function, $Z$ agrees with a $C$-fundamental vector field $V$.
Let $\phi$ be a defining function for the zero locus of $V$.
We may then assume that $Z=\frac{1}{\phi}V$.
Then
\[
fd\eta_i(Y,Z)=fYd\eta_i\frac{1}{\phi}V
=- f\frac{Y\phi}{\phi^2}(\eta_iV) + \frac{f}{\phi}Y(\eta_iV)
\]
Since $\eta_i$ is dual to $V_i$, the function $\eta_iV$ is constant
and the second term vanishes.
Therefore
\[
fd\eta_i(Y,Z)
=- \left(\frac{f}{\phi}\right)
\left(\frac{Y\phi}{\phi}\right)\left(\eta_iV\right).
\]
We claim that all three factors are holomorphic.
Indeed, the zero locus of $\phi$ is contained in $D$ 
(with multiplicity one) and
$f$ vanishes on $D$. By construction the horizontal vector field
$V$ is tangent to $D$, hence $V(\phi)$ vanishes along the zero locus
of $\phi$ (which is contained in $D$) and therefore $V\phi/\phi$
is holomorphic. Finally $\eta_iV$ is evidently holomorphic, since it
is a constant function.
\end{proof}
\begin{proof}[Proof of theorem \ref{thm3}]
The obvious connection on the trivial bundle $H\times C\to H$
induces a flat connection on $H\times(\bar C/\sim)\to Y=H/\Lambda$.

Hence the statement follows from the preceding theorem (thm.~\ref{thm2}).
\end{proof}

\end{document}